\documentclass[a4paper,12pt]{article}
\usepackage{amssymb}
\usepackage{amsfonts}
\usepackage{amsmath}
\usepackage{color}
\usepackage{verbatim}
\usepackage{latexsym,array,theorem,bbm}
\usepackage[colorlinks=true,linkcolor=blue,citecolor=magenta]{hyperref}    

\newcommand{\RR}{\mathbb{R}}

\newcommand{\NN}{\mathbb{N}}
\newcommand{\QQ}{\mathbb{Q}}

\newcommand{\cB}{{\mathcal B}}
\newcommand{\cC}{{\mathcal C}}

\renewcommand{\P}{\mathsf{P}}
\newcommand{\cov}{\mathsf{Cov}}
\newcommand{\cum}[1][]{\mathsf{Cum}_{\,#1}}
\newcommand{\EE}{\mathsf{E}}
\newcommand{\PP}{\mathsf{P}}
\newcommand{\Dn}{\mathsf{D}^2}
\newcommand{\D}{\mathsf{D}}
\newcommand{\proofend}{\hfill\mbox{$\Box$}\medskip}

\newcommand{\Hexp}{\Gamma}

\numberwithin{equation}{section}

\theoremstyle{change}
\newtheorem{Lem}{Lemma.}[section]
\newtheorem{Thm}[Lem]{Theorem.}

\newtheorem{Cor}[Lem]{Corollary.}
\newtheorem{Def}[Lem]{Definition.}
\newtheorem{Ass}[Lem]{Assumption.}
\newcounter{Not}

\newtheorem{Rem}[Lem]{Remark.}
\newtheorem{Ex}[Lem]{Example.}

\begin{document}

\frenchspacing

\title{\bfseries\Large Path properties of dilatively stable processes\\ and singularity of their distributions}
\author{{\sc\large Endre $\text{Igl\'oi}^{\ast,\diamond}$}
         {\large and}
         {\sc\large M\'aty\'as $\text{Barczy}^{*}$}  }
\date{}
\maketitle

\vskip0.2cm
* University of Debrecen, Faculty of Informatics, Pf.~12, H--4010 Debrecen, Hungary,
Tel.: 06-52-512900, Fax: 06-52-512996

e-mails:  igloi@tigris.unideb.hu (E. Igl\'oi), barczy.matyas@inf.unideb.hu (M. Barczy),

$\diamond$ Corresponding author.

\vskip0.2cm




\vspace*{5mm}

\date{}

\pagestyle{myheadings}

\markboth{Dilatively stable processes: path properties and singularity}
         {Dilatively stable processes: path properties and singularity}

\begin{abstract}
First, we present some results about the H\"older continuity of the sample paths
of so called dilatively stable processes
which are certain infinitely divisible processes having a more general scaling property
than self-similarity.
As a corollary, we obtain that
 the most important \!$(H,\delta)$\!-dilatively stable limit processes
(e.g., the LISOU and the LISCBI processes,  see Igl\'oi \cite{Igl})
 almost surely have a local H\"older exponent \!$H.$\!
 Next we prove that, under some slight regularity assumptions,
 any two dilatively stable processes with stationary increments
 are singular (in the sense that their distributions have disjoint supports)
if their parameters \!$H$\! are different.
We also study the more general case of not having stationary increments.
Throughout the paper we specialize our results to some basic dilatively stable
 processes such as the above-mentioned limit processes and the fractional L\'evy process.

\end{abstract}

\section{Introduction}
\label{Isect}

Path properties are important features of a stochastic process, this fact does not call for an explanation.
Given a set of stochastic processes, the same is true for the singularity
of their distributions on the space of sample paths
(by which we understand that under certain conditions any two such distributions have disjoint supports).

A well-known feature enabling a process to have nice path properties is self-similarity.
Systematic investigations of the sample path properties of self-similar processes had
 been begun by Vervaat's classical paper \cite{Ver},
 and continued by Vervaat \cite{Ver1},  Takashima \cite{Tak} and Watanabe and Yamamuro
 \cite{WatYam} (the latter one gives laws of the  iterated logarithm
 for multidimensional self-similar processes with independent increments).
 Maejima \cite{Maejima1983},  Kono and Maejima \cite{KonMae} and Samorodnisky \cite{Sam}
 (among others) treated the sample path properties of some special self-similar processes.
However, there are very few results in the literature about the singularity properties
 of distributions of self-similar processes.
To the authors' knowledge the only paper concerning this subject is that of Prakasa Rao \cite{Rao},
 which states that the distributions of two fractional Brownian motions (FBMs) with different
 Hurst parameters are singular with respect to each other.
The question of singularity for other types of processes, for instance, for diffusion processes
 deserved more attention, see, e.g.,
 Ben-Ari and Pinsky \cite{BenPin}, Jacod and Shiryaev \cite{JacShi} and Zhang \cite{Zha}.

In this paper we show that the so-called dilatively stable processes, introduced by Igl\'oi \cite{Igl},
 can have nice sample path and singularity properties.
Dilative stability is a generalization of self-similarity. For the comparison of these
 two notions, we give both definitions.

\begin{Def}\label{DefSS}
Let $\alpha>0.$  A process $\{X(t),\,t\geq0\}$
starting from zero (i.e., $X(0)=0$\!)
is called $\alpha$\!-self-similar if it is not identically zero and fulfills the scaling relation
\begin{align}\label{ssrel}
\forall\, T>0 : \, X(T\,\cdot\,)\overset{\text{\text{fd}}}\sim T^\alpha X(\,\cdot\,),
\end{align}
where $\overset{\text{\text{fd}}}{\sim }$ denotes that the finite-dimensional distributions are the same.
\end{Def}

Dilative stability, see Igl\'oi \cite{Igl}, is an analogous property of certain infinitely divisible processes
(all finite-dimensional distributions are infinitely divisible)
 involving a scaling also in the convolution exponent.
The exact definition is as follows.

\begin{Def}\label{DefDS}
Let $\alpha>0$ and $\delta\leq2\alpha.$
A process $\{X(t),\,t\geq0\}$
starting from zero
is said to be $(\alpha,\delta)$\!-dilatively stable
if its finite-dimensional distributions are infinitely divisible,
$X(1)$ is non-Gaussian, $X(t)$ has finite moments of all orders for all $t\geq 0$,
and fulfills the scaling relation
\begin{equation}
\label{dsrel}
\forall\, T>0: \,
 X(T\,\cdot\,)\overset{\text{fd}}\sim T^{\alpha-\delta/2} X^{\circledast T^\delta}(\,\cdot\,).
\end{equation}
Here for all $c>0$, we denote by
 $X^{\circledast c}=\{X(t),\,t\geq0\}^{\circledast c}$ the $c$\!-th convolution power of
 $\{X(t),\,t\geq0\},$ that is, $\{X(t),\,t\geq0\}^{\circledast c}$ is a process
 the finite dimensional distributions  of which are the $c$\!-th convolution powers of the corresponding ones of
 $\{X(t),\,t\geq0\}.$
\end{Def}

In Appendix \ref{appendix} we give more insight to the properties of a dilatively stable process that
 are supposed in Definition \ref{DefDS}.
We also note that in Kaj \cite[Section 3.6]{Kaj} one can find a somewhat similar,
 but not so general concept called aggregate-similarity.
One of the main differences between the concept of dilative stability and aggregate-similarity
 is that the later one defined only in the case of $n$\!-th convolution powers with $n\in\NN$
 and when the relation $2\alpha-\delta=2$ holds for the parameters
 (which is in fact the most important case, see the parametrization of the LIS processes in Example \ref{LISexample}).

Let us observe that  a not identically zero process $\{X(t),\,t\geq0\}$ starting
 from zero and fulfilling the scaling relation \eqref{dsrel} with $\alpha>0$ and
 $\delta=0$ is an $\alpha$-self-similar process, roughly speaking $(\alpha,0)$\!-dilative stability
 is just \!$\alpha$\!-self-similarity.
Hence dilative stability is a generalization of self-similiraty.
Many results, such as Lamperti's theorems, can be transferred from the self-similar
to the dilatively stable case, see Igl\'oi \cite{Igl}.
Self-similar processes are  fixed points of their renormalization operators
(see Taqqu \cite{Taq}),
and the same is true for dilatively stable processes  (see, Igl\'oi \cite[Theorem 2.8.4 (DS)]{Igl}).
Accordingly, the possible limit processes in linear  rescaling (i.e.,
 in self-similar renormalization limit theorems) are the self-similar processes,
  see Lamperti \cite{Lam},
 and if there is a rescaling also in the convolution exponent (i.e., in dilatively stable
 renormalization limit theorems) the possible limit processes are the dilatively stable ones,
 see, Igl\'oi \cite[Theorem 2.2.7]{Igl}.
This is why dilative stability is important.

Important examples for dilatively stable processes are non-Gaussian {\bf f}ractional {\bf L}\'evy {\bf p}rocesses
 (FLPs) (i.e., FLPs, where the underlying two-sided L\'evy process is non-Gaussian, but possibly with a Gaussian component)
 having zero mean and finite moments of all orders.
FLPs were originally introduced by Benassi et al. \cite{Ben} and Marquardt \cite{Mar}.
 For historical fidelity we note that in \cite{Ben} a FLP is called a moving-average fractional L\'evy motion.
We also recall that Marquardt \cite[Theorem 4.4]{Mar} proved that a FLP with an underlying two-sided L\'evy process having zero mean,
 finite second moment and not having a Brownian component cannot be self-similar.
However, by Igl\'oi \cite[Example 2.1.7]{Igl}, a non-Gaussian FLP (considered only on $[0,\infty)$\!) having zero mean
 and finite moments of all orders is $(H,1)$-dilatively stable with stationary increments, where $H\in(1/2,1)$ is
 the so-called Hurst parameter or long memory parameter (see also Kaj \cite[page 212]{Kaj}).
Roughly speaking, non-Gaussian FLPs are not self-similar but they belong to a wider class of processes, to the class of dilatively stable
 processes, which also underlines the importance of dilative stability.

The main contributions of this paper relate to dilatively stable processes with
 stationary increments, the general properties of which are treated by Igl\'oi \cite[Section 2.7]{Igl}.
Here we point out that in this case the range of the parameter \!$H\doteq\alpha$\! is
 $ (0,1],$  see, Igl\'oi \cite[Theorem 2.7.1 (DS) 1)]{Igl}.
In case $H=1,$ i.e., in case of a $(1,\delta)$\!-dilatively stable process we have
 $\delta=0$ (i.e., the process is self-similar) and it takes the form
 $X(t) = tX(1),$ $t\geq 0,$ almost surely (a.s.), i.e., it is degenerate.
Indeed, by Igl\'oi \cite[Theorem 2.7.1 (DS) 2)]{Igl}, the process is degenerate, and
 if we suppose on the contrary that $\delta\ne 0,$ then, by \eqref{dsrel}
 (with $T\doteq c^{1/\delta},$ $c>0$),
 \[
    X(c^{1/\delta})\sim c^{\frac{1}{\delta}(1-\delta/2)} X^{\circledast c}(1),
    \qquad c>0,
 \]
 which yields that
 \[
    c^{\frac1\delta} X(1) \sim c^{\frac{1}{\delta}- \frac{1}{2}} X^{\circledast c}(1),
    \quad \text{i.e.,}\quad
    X(1) \sim \frac{1}{\sqrt{c}} X^{\circledast c}(1),
    \qquad  c>0.
 \]
Hence $X(1)$ is Gaussian with mean zero (which follows by considering characteristic functions),
in contradiction with Definition \ref{DefDS}.

Moreover, an $(H,\delta)$\!-dilatively stable process  $\{X(t),\,t\geq 0\}$ with
 stationary increments has the same covariance function as a FBM with parameter \!$H$\!
 (apart from a constant factor):
 \[
  \cov(X(t_1),X(t_2)) = \frac{1}{2}\D^2X(1)\left(t_1^{2H} + t_2^{2H} - \vert t_1-t_2\vert^{2H}\right),\qquad t_1,t_2\geq 0,
 \]
 see Igl\'oi \cite[Theorem 2.7.2]{Igl} (where the parameter \!$H$\! can not be one, however the proof given
 there works also in case $H=1$\!).

Moments and cumulants provide a very handy tool throughout when dealing with dilative stability.
Using this tool in  Section 2 we obtain some results for the almost sure limit behavior
 of dilatively stable processes at zero, at infinity and, in case of processes with stationary increments,
 at any point, see Lemmas \ref{Lem1}, \ref{Leminfty} and \ref{Lem2}.
Then, by the help of the Kolmogorov--Chentsov theorem,
we characterize dilatively stable processes with stationary increments
from the point of view of H\"older continuity of their sample paths,  see Theorem \ref{Holder-thm}.
This characterization has a well-known self-similar analogue presented in Remark \ref{REM_ss1}.
Section 3 contains  some results about the singularity of the distributions of dilatively
 stable processes, induced on the space of continuous functions.
Though a dilatively stable process does not admit automatically
 continuous sample paths almost surely, in the last section we will restrict ourselves to such processes.
We note that an $(H,\delta)$\!-dilatively stable process with stationary increments and zero
 mean such that either $\delta\ne 2H$ or $H>1/2$ has a continuous modification,
 see Corollary \ref{C-beliseg}.
Using the path properties proved in Section 2, we will find that, under some slight regularity
 conditions, the distributions of dilatively stable processes  with
 different parameters $\alpha$ are pairwise singular, see
 Theorem \ref{Thm1} and Theorem \ref{Thm2} (this latter one is for the stationary increments case).
 Theorem \ref{Thm2} can be carried out to self-similar processes with finite  absolute moment
 and stationary increments, implying, particularly, a new (and simple) proof of the pairwise
 singularity of FBMs with different Hurst parameters, see Theorem \ref{Thm2ss}.

Throughout the paper we will specialize our results to some particular dilatively stable
 processes with stationary increments presented in the next example.

\begin{Ex}\label{LISexample}
The following processes are dilatively stable with stationary increments, see Igl\'oi \cite{Igl}:\newline
$\bullet$ LISOU process ({\bf l}imit of {\bf i}ntegrated {\bf s}uperposition of {\bf O}rnstein--{\bf U}hlenbeck type processes),\newline
$\bullet$ LISCBI process ({\bf l}imit of {\bf i}ntegrated {\bf s}uperposition of {\bf c}ontinuous state
{\bf b}ranching processeses with {\bf i}mmigration),\newline
$\bullet$ LISDLG process ({\bf l}imit of {\bf i}ntegrated {\bf s}uperposition of {\bf d}iffusion processes with
{\bf l}inear {\bf g}enerator; this is a particular LISCBI process),\newline
$\bullet$  non-Gaussian FLP ({\bf f}ractional {\bf L}\'evy {\bf p}rocess), i.e.,
 FLP defined on $[0,\infty),$ where the underlying two-sided L\'evy process is non-Gaussian (but possibly with a Gaussian
 component) having zero mean and finite moments of all  orders.
\newline
\indent These processes have a parameter \!$H,$\! which can now be only
 in the interval $(1/2,1),$ hence it can be called the Hurst parameter or long memory parameter.
The LISOU, LISCBI and LISDLG processes with parameter \!$H$\! are  $(H,2H{\,-\,}2)$\!-dilatively stable,
 while the non-Gaussian FLP with parameter \!$H$\! is $(H,1)$\!-dilatively stable,
and all of these processes have the same covariance function as a FBM with  parameter $H.$
\end{Ex}

\medskip

We will suppose throughout that the processes have zero mean.
This is a natural assumption when dealing with path properties, since the mean function,
 as a deterministic function, can be handled separately.
One can also check that the subtraction of the mean function preserves dilative stability.
In addition, if an $(H,\delta)$\!-dilatively stable process has stationary increments,
then its mean function is automatically zero,
unless $H+\delta/2=1,$ see Igl\'oi \cite[Theorem 2.7.1 (DS) 3)]{Igl}.

\section{Path properties}
\label{pathsect}

Throughout in this paper $I\subseteq\mathbb[0,\infty)$ denotes an interval.
Let $\gamma\in(0,1].$
A function $f:I\rightarrow\mathbb{R}$ is called \emph{locally \!$\gamma$\!-H\"older continuous}
if for every bounded subinterval $J\subseteq I,$
\[
\sup_{{t,s\in J,}\atop{t\neq s}}\frac{\vert f(t)-f(s)\vert}{\vert t-s\vert^\gamma}<\infty,
\]
see, e.g., {Revuz and Yor \cite[page 26]{RevYor}.
Clearly, local \!$\gamma_1$\!-H\"older continuity implies local \!$\gamma_2$\!-H\"older continuity
 if  $1\geq \gamma_1\geq\gamma_2>0.$
This relation gives rise to the following notion of local H\"older exponent.
The value
\[
\Hexp_f\doteq\sup_{\gamma\in(0,1]}\big\{\gamma: f
                           \text{\ is locally \!$\gamma$\!-H\"older continuous}\big\}
\]
is called the \emph{(optimal) local H\"older exponent} of  a function \!$f$
 and set to $0$ if $f$ is not locally H\"older continuous.
This notion is similar to the notion of the optimal H\"older index at a point, see, e.g.,
 Jaffard \cite{Jaf} or Fleischmann, Mytnik and Wachtel \cite{FleMytWac}.}

 In what follows by the expression that `an infinitely divisible  distribution has
 a Gaussian component' we mean that in its L\'evy-Khintchine representation the Gaussian part
 has positive variance.

In some cases we will refer to the following assumptions.

\begin{Ass}
\label{Assum_delta+}
If $\{X(t),\,t\geq0\}$ is an $(\alpha,\delta)$\!-dilatively stable process with $\delta\geq 0,$
 then the distribution of $X(1)$ (equivalently, the distribution of $X(t)$ for any $t>0$\!)
 has a Gaussian component.
\end{Ass}

\begin{Ass}
\label{Assum_delta-}
If $\{X(t),\,t\geq0\}$ is an $(\alpha,\delta)$\!-dilatively stable process with $\delta\leq 0,$
 then the distribution of $X(1)$ (equivalently, the distribution of $X(t)$ for any $t>0$\!)
 has a Gaussian component.
\end{Ass}

Note that the above two assumptions are the dual of each other with respect to \!$\delta,$
and we will always indicate explicitly which of them is used.
Furthermore, when Assumption \ref{Assum_delta+} is supposed and the parameter
 \!$\delta$\! of an $(\alpha,\delta)$\!-dilatively stable process
 is negative, then this assumption does not come into play,
 and a similar statement holds for Assumption \ref{Assum_delta-}
 and an $(\alpha,\delta)$\!-dilatively stable process with positive \!$\delta.$
The intrinsic reason for considering only the non-negative values $\delta$
 in Assumption \ref{Assum_delta+} is that the LISOU, LISCBI and LISDLG processes
 (the dilatively stable processes with stationary increments presented in
 Example \ref{LISexample}) may not have a Gaussian component and their parameter
 $\delta$ is negative.
 However, this also means that Assumption \ref{Assum_delta-} may not hold for these particular
 dilatively stable processes, and from this specific point of view Assumptions \ref{Assum_delta+}
 and \ref{Assum_delta-} are not exactly the dual of each other.
On the other side, we call the attention that these assumptions are not too restrictive,
 since, as it is easy to see, the independent sum of a Gaussian $\alpha$\!\!-self-similar process
 with zero mean and an $(\alpha,\delta)$\!-dilatively stable process remains
 $(\alpha,\delta)$\!-dilatively stable.
Particularly, the independent sum of an $(H,\delta)$\!-dilatively stable process and a FBM with parameter
 \!$H$\! remains $(H,\delta)$\!-dilatively stable, see  Igl\'oi \cite[Proposition 2.7.5]{Igl}.
Thus, we can make a Gaussian component preserving the dilative stability.
Finally,  we remark that under the non-trivial cases of Assumption \ref{Assum_delta+} or
 Assumption \ref{Assum_delta-}, i.e., when there exists a Gaussian component,
 the distribution of $X(1)$ is absolutely continuous  (see, e.g., Sato \cite[Lemma 27.1]{Sat}),
 and hence $\P(X(1)=0)=0$ in this case.

The following lemma treats the sample path behaviour at zero of a dilatively stable process.
In what follows by a zero-sequence we mean a sequence of real numbers converging to 0.

\begin{Lem}\label{Lem1}
There exists a zero-sequence $(t_n)_{n\in\mathbb{N}}$ with positive terms such that
for any $(\alpha,\delta)$\!-dilatively stable process $\{X(t),\,t\geq0\}$ with zero mean,
the following assertions hold.\newline
(i) If $\kappa<\alpha,$ then
\begin{equation}
\label{Halatt}
\limsup_{n\to\infty}\frac{\vert X(t_n)\vert}{t_n^\kappa}=0\quad\text{a.s.}
\end{equation}
(ii) If $\kappa>\alpha$ and  Assumption \ref{Assum_delta+} holds, then
\begin{equation}
\label{Hfolott}
\limsup_{n\to\infty}\frac{\vert X(t_n)\vert}{t_n^\kappa}=\infty\quad\text{a.s.}
\end{equation}
Further, for any sequence $(t_n)_{n\in\mathbb{N}}$ with positive terms such that
\begin{equation}\label{limsupfelt}
\limsup_{n\to\infty}\sqrt[n]{t_n}<1,
\end{equation}
the assertions of parts (i) and (ii) hold.
\end{Lem}

\noindent{\bf Proof.}
Let $(t_n)_{n\in\mathbb{N}}$ be a sequence with positive terms such that \eqref{limsupfelt} holds.
Then the root test yields that the series  $\sum_{n=1}^\infty t_n$  is convergent,
and hence $(t_n)_{n\in\mathbb{N}}$ is a zero-sequence.

(i) By (\ref{dsrel}), we have
\[
\Dn X(t_n)=\Dn\big(t_n^{\alpha-\delta/2}X^{\circledast t_n^\delta}(1)\big)
=t_n^{2\alpha-\delta+\delta}\Dn X(1)
=t_n^{2\alpha}\Dn X(1), \quad  n\in\NN.
\]
Hence the Markov inequality yields that
\[
\sum_{n=1}^\infty \P\left(\frac{\vert X(t_n)\vert}{t_n^\kappa}>\varepsilon\right)
\leq \frac{\Dn X(1)}{\varepsilon^2}\, \sum_{n=1}^\infty t_n^{2(\alpha-\kappa)}<\infty,
 \quad   \varepsilon>0,
\]
 where the convergence of the series is a consequence of (\ref{limsupfelt}).
Applying the Borel--Cantelli lemma we obtain that $\lim_{n\to\infty}X(t_n)/t_n^\kappa=0$ a.s.,
 which implies (\ref{Halatt}).
We remark that in this case any sequence $(t_n)_{n\in\NN}$ with positive terms
 satisfying \eqref{limsupfelt} is obviously universal.

(ii) Using (\ref{dsrel}), we have
\begin{equation}
\label{negyedik}
\frac{\vert X(t_n)\vert}{t_n^\kappa}
\sim t_n^{\alpha-\kappa}\,\frac{\left\vert X^{\circledast t_n^\delta}(1)\right\vert}{t_n^{\delta/2}},
  \quad  n\in\NN.
\end{equation}
First we show that the right-hand side (and hence the left-hand side)  of \eqref{negyedik} converges
 in probability to infinity  along some appropriate subsequence.
We consider the three cases corresponding to the sign of $\delta$ separately.

\noindent$\bullet$ If $\delta<0,$ then $X^{\circledast t_n^\delta}(1)\big/t_n^{\delta/2}$ converges
 in distribution to a normal distribution with variance $\EE (X(1))^2.$
Indeed, $\lim_{n\to\infty}t_n^\delta=\infty,$ and
 \[
      X^{\circledast t_n^\delta}(1)
       \sim \xi_1+\cdots+\xi_{\lfloor t_n^\delta\rfloor}+\eta_n,
 \]
 where  $\lfloor\cdot\rfloor$ denotes the integer part,
 $\xi_1,\ldots,\xi_{\lfloor t_n^\delta\rfloor}$ and $\eta_n$  are independent random variables
 such that  $\xi_1,\ldots,\xi_n$  has a common distribution as $X(1)$ has and
 $\eta_n\sim X^{\circledast (t_n^\delta -\lfloor t_n^\delta\rfloor)}(1).$
By the central limit theorem we get that
 \[
      \frac{\xi_1+\cdots+\xi_{\lfloor t_n^\delta\rfloor}}{t_n^{\delta/2}}
 \]
 converges in distribution to a normal distribution with variance $\EE (X(1))^2,$
 while the remainder term
 $\eta_n/t_n^{\delta/2}$ converges to 0 in $L^2$ as $n\to\infty.$
Therefore the right-hand side of (\ref{negyedik}) converges in probability to infinity,
 and hence so does the left-hand side.
(Here we used the fact that
 if $\zeta_n,\, n\in\NN,$ and \!$\zeta$\! are non-negative random variables
 such that \!$\zeta_n$\! converges to \!$\zeta$\! in distribution,
 $\zeta$ has a continuous distribution function  [hence $\PP(\zeta>0)=1$\!] and
 $(c_n)_{n\in\NN}$ is a sequence of positive real numbers converging
 to infinity, then $c_n\zeta_n$ converges in probability to infinity as $n\to\infty$,
 i.e., $\forall\,M>0:\,\lim_{n\to\infty}\PP(c_n\zeta_n<M)=0.$\!)

\noindent $\bullet$ If $\delta=0,$ then the right-hand side of (\ref{negyedik}) is
 $t_n^{\alpha-\kappa}\vert X(1)\vert.$
Using Assumption \ref{Assum_delta+} and the fact that if
 at least one term of an independent  (finite) sum  of random variables
 has a continuous distribution, then also the sum itself has a continuous distribution
 (see, e.g., Sato \cite[Lemma 27.1]{Sat}),
 we have $\PP(X(1)=0)=0,$ i.e., $X(1)$ has no atom at zero.
 Hence $t_n^{\alpha-\kappa}\vert X(1)\vert$  converges almost surely to infinity,
 which yields that the left-hand side of (\ref{negyedik}) converges in probability to $\infty.$

\noindent $\bullet$ If $\delta>0,$ then, since
\begin{equation}\label{marad}
 \Dn\left(\frac{X^{\circledast t_n^\delta}(1)}{t_n^{\delta/2}}\right)
   =\Dn X(1)
   =\EE(X(1))^2,\quad n\in\mathbb{N},
\end{equation}
the tightness of the sequence of distributions of
$X^{\circledast t_n^\delta}(1)/t_n^{\delta/2},$ $n\in\mathbb{N},$ follows by the Markov inequality.
 Indeed, with the notation
 \[
   K_\varepsilon
      \doteq\Big[ -\!\sqrt{\EE(X(1))^2/\varepsilon},\ \sqrt{\EE(X(1))^2/\varepsilon}\,\Big],
      \quad \varepsilon>0,
 \]
 for all $n\in\NN,$ we have
 \[
    \PP\left( \frac{X^{\circledast t_n^\delta}(1)}{t_n^{\delta/2}}
              \in \RR\setminus K_\varepsilon \right)
     = \PP \left( \left\vert\frac{X^{\circledast t_n^\delta}(1)}{t_n^{\delta/2}}\right\vert
                  >  \sqrt{\frac{\EE(X(1))^2}{\varepsilon}} \right)
     \leq \varepsilon.
 \]
Therefore, by Prohorov's theorem, there exists a subsequence
 $X^{\circledast t_{n_k}^\delta}(1)/t_{n_k}^{\delta/2},$ $k\in\mathbb{N},$ which converges in distribution.
Using Assumption \ref{Assum_delta+} and the L\'evy--Khintchine formula we get that for all $k\in\NN,$
 $X^{\circledast t_{n_k}^\delta}(1)/t_{n_k}^{\delta/2}$ has a Gaussian component,
 the distribution of which is the same as
 the distribution of the Gaussian component of $X(1).$
Indeed, if \!$G$\! denotes the Gaussian component of $X(1),$ then
 $G^{\circledast t_{n_k}^\delta}/t_{n_k}^{\delta/2}$ is the Gaussian component of
 $X^{\circledast t_{n_k}^\delta}(1)/t_{n_k}^{\delta/2}$ and it has the same distribution as \!$G$\!
 (here we also use that \!$G$\! has zero mean since $X(1)$ has zero mean).
Thus the limiting distribution of this subsequence also has a Gaussian component.
Therefore the limiting distribution is continuous and hence it
 has no atom at zero, and, as it was explained earlier in the proof of the case $\delta<0,$
 this fact ensures that the right-hand side of
 (\ref{negyedik}) converges along the subsequence $(t_{n_k})_{k\in\mathbb{N}}$ in probability
 to infinity, hence so does the left-hand side (along the subsequence $(t_{n_k})_{k\in\mathbb{N}}$).

 We call the attention that in the proof of the case $\delta>0$
(i.e., in the last case) we do not use that \!$\delta$\!
 is positive, but Assumption \ref{Assum_delta+} comes into play in this case.
This also shows that the proof of the case $\delta>0$ remains still valid for the case $\delta<0$
 under the additional assumption that $X(1)$ has a Gaussian component.
 (However, the proof of the case $\delta<0$ presented above does not work for the case $\delta>0.$\!)

In each of the above three cases (denoted by bullets) we obtained that
a subsequence of the left-hand side of (\ref{negyedik}) (in case of $\delta\leq 0$ the whole
 sequence) converges in probability to infinity.
Thus, by the Riesz lemma, there is some subsequence of the left-hand side of (\ref{negyedik}),
 which converges to infinity almost surely.
Indeed, Riesz's lemma can be proved for a sequence of random variables $(\zeta_n)_{n\in\NN}$
 converging in probability to $\infty,$ as follows.
Since for all $k\in\NN$ there exists some $n_k\in\NN$ such that $\P(\zeta_{n_k}<k)<1/2^k,$
 the Borel--Cantelli lemma yields that
 \[
    \P(\zeta_{n_k}<k \;\; \text{for infinitely many} \;\; k\geq 1)=0.
 \]
Hence $\P(\liminf_{k\to\infty}\zeta_{n_k}=\infty)=1$ which yields that
 $\P(\lim_{k\to\infty}\zeta_{n_k}=\infty)=1$.
This implies (\ref{Hfolott}) and we also have the universality of any sequence $(t_n)_{n\in\NN}$
 with positive terms satisfying \eqref{limsupfelt}.
\proofend

\begin{Rem}
We note that in part (i) of Lemma \ref{Lem1} one can also write $\lim$ instead of $\limsup.$
The reason for writing $\limsup$ is that we will use only this later on.
We also remark that the case $\kappa=\alpha$ is not covered by Lemma \ref{Lem1},
 since we do not need it and we can not address any result  in this case.
\end{Rem}

The following lemma shows that the behaviour of the sample paths at infinity can be characterized similarly
 as their behaviour at zero.

\begin{Lem}\label{Leminfty}
There exists a sequence $(t_n)_{n\in\mathbb N}$ converging to infinity such that for any
$(\alpha,\delta)$\!-dilatively stable process  $\{X(t),\,t\geq0\}$  with zero mean the following assertions
hold.
\newline
(i) If $\kappa<\alpha$ and Assumption \ref{Assum_delta-} holds, then
\begin{equation}
\label{Halatt1}
\limsup_{n\to\infty}\frac{\vert X(t_n)\vert}{t_n^\kappa}=\infty\quad\text{a.s.}
\end{equation}
(ii) If $\kappa>\alpha,$  then
\begin{equation}
\label{Hfolott1}
\limsup_{n\to\infty}\frac{\vert X(t_n)\vert}{t_n^\kappa}=0\quad\text{a.s.}
\end{equation}
Further, for any sequence $(t_n)_{n\in\mathbb{N}}$ with positive terms such that
 $\limsup_{n\to\infty}\!\sqrt[n]{1/t_n}<1,$ the assertions of parts (i) and (ii) hold.
\end{Lem}

\noindent{\bf Proof.} One can use the same arguments as in the proof of Lemma \ref{Lem1},
 but with some changes.
Namely, the sequence $(t_n)_{n\in\mathbb N}$ can be the reciprocal of the sequence
(with positive terms satisfying \eqref{limsupfelt}) used in the proof of Lemma \ref{Lem1}.
Then the proof of statement (ii) corresponds to that of statement (i) in Lemma \ref{Lem1},
while part (i) goes the same way as part (ii) in Lemma \ref{Lem1},
just the  cases $\delta<0$ and $\delta>0$ have to be interchanged.
\proofend

\medskip

In the stationary increments case Lemma \ref{Lem1} can be formulated not only for zero-sequences
 but also for sequences converging to a non-negative real number, see as follows.

\begin{Lem}
\label{Lem2}
For any $t_0\geq0$ there exists a sequence $(t_n)_{n\in\mathbb{N}}$ converging to $t_0$ such that
$t_n\neq t_0,$ $n\in\mathbb{N},$ and for any $(H,\delta)$\!-dilatively stable process $\{X(t),\,t\geq0\}$ with stationary increments
and zero mean the following assertions hold.\newline
(i) If $\kappa<H,$ then
\begin{equation}
\label{Halatts}
\limsup_{n\to\infty}\frac{\vert X(t_n)-X(t_0)\vert}{\vert t_n-t_0\vert^\kappa}=0\quad\text{a.s.}
\end{equation}
(ii) If $\kappa>H$ and Assumption \ref{Assum_delta+} holds, then
\begin{equation}
\label{Hfolotts}
\limsup_{n\to\infty}\frac{\vert X(t_n)-X(t_0)\vert}{\vert t_n-t_0\vert^\kappa}=\infty\quad\text{a.s.}
\end{equation}
Further, for any sequence $(t_n)_{n\in\mathbb{N}}\doteq(t_0+\tilde t_n)_{n\in\mathbb{N}},$
 where $(\tilde t_n)_{n\in\mathbb{N}}$ is a sequence with positive terms such that
 $\limsup_{n\to\infty}\sqrt[n]{\tilde t_n}<1,$ the assertions of parts (i) and (ii) hold.
\end{Lem}

\noindent{\bf Proof.}
The process $Y(t)\doteq X(t+t_0)-X(t_0),$ $t\geq0,$ is $(H,\delta)$\!-dilatively stable  with zero mean
 such that the distribution of $Y(1)$ has a Gaussian component
 if $X(1)$ has.
Hence Lemma \ref{Lem1} yields that there exists a zero-sequence $(\tilde t_n)_{n\in\mathbb{N}}$
 with positive terms such that if $\kappa<H,$ then
 \[
   \limsup_{n\to\infty}\frac{\vert X(\tilde t_n + t_0)-X(t_0)\vert}
                          {\tilde t_n^\kappa}=0\quad\text{a.s.,}
 \]
and if $\kappa>H$ and Assumption \ref{Assum_delta+} holds, then
\begin{align*}
 \limsup_{n\to\infty}\frac{\vert X(\tilde t_n + t_0)-X(t_0)\vert}
                          {\tilde t_n ^\kappa}=\infty\quad\text{a.s.}
\end{align*}
With the definition $t_n\doteq\tilde t_n + t_0,$ $n\in\NN,$ we have the assertions
 of the lemma.
(This also shows that in \eqref{Halatts} and \eqref{Hfolotts} one can write
 $(t_n-t_0)^\kappa$ instead of $\vert t_n-t_0\vert^\kappa.$\!)
\proofend

\medskip

The above lemmas, interesting in their own rights, will be used mainly in the next section,
 but Lemma \ref{Lem2} will appear also in the proof of the following theorem.
The rest of the section deals with H\"older continuity of the sample paths of
 dilatively stable processes with stationary increments.

\begin{Thm}\label{Holder-thm}
Let $\{X(t),\,t\geq0\}$ be an $(H,\delta)$\!-dilatively stable process with stationary increments and zero mean
 such that either $\delta\neq2H$ or $H>1/2.$
Then $\{X(t),\,t\geq0\}$ has a continuous modification, the sample paths of which
are locally \!$\gamma$\!-H\"older continuous
\renewcommand{\labelenumi}{{\rm(\roman{enumi})}}
\begin{enumerate}
 \item for every $\gamma\in(0,H)$ if $\delta<0,$
 \item for every $\gamma\in(0,H-\delta/2)$ if $0\leq\delta<2H,$
 \item for every $\gamma\in(0,H-1/2)$ if $\delta=2H$ and $H>1/2.$
\end{enumerate}
Moreover, under Assumption \ref{Assum_delta+}, for the local H\"older exponent  $\Hexp_X$
 of the sample paths of the above continuous modification of the process $\{X(t),\,t\geq0\}$ we have
\begin{enumerate}
 \item[$(\widetilde i)$] $\Hexp_X=H$ if $\delta<0,$
 \item[$(\widetilde {ii})$] $\Hexp_X\in[H-\delta/2, H]$ if $0\leq\delta<2H,$
 \item[$(\widetilde {iii})$] $\Hexp_X\in[H-1/2,H]$ if $\delta=2H$ and $H>1/2.$
\end{enumerate}
\end{Thm}

\noindent{\bf Proof.}
We are going to apply the Kolmogorov--Chentsov theorem see, e.g.,
 Revuz and Yor \cite[Chapter I, Theorem 2.1]{RevYor}.
Therefore we have to prove that $\{X(t),\,t\geq0\}$ satisfies Kolmogorov's condition:
 there exist constants $c,p,q>0$ such that
\begin{equation}
\label{Kolmo}
 \EE\big\vert X(t)-X(s)\vert^p\leq c\vert t-s\vert^{1+q}, \qquad s,t\geq 0.
\end{equation}
Then, by the Kolmogorov--Chentsov theorem,
$\{X(t),\,t\geq0\}$  has a continuous modification which  is locally
\!$\gamma$\!-H\"older
continuous for every $\gamma\in(0,q/p).$
Clearly, it is sufficient for (\ref{Kolmo}) to hold for every
$s,t\geq0$ for which $0\leq t-s<1,$ so in what follows we suppose that \!$t$\! and \!$s$\! are
of these kinds.
Indeed, in this case for all $n\in\NN\cup\{0\},$
 \[
    \EE\big\vert X^{(n)}(t)-X^{(n)}(s)\vert^p\leq c\vert t-s\vert^{1+q}, \qquad s,t\geq 0,
 \]
 where
 \[
    X^{(n)}(t)
      \doteq\begin{cases}
          X(n/2) & \text{if \ $t\leq n/2,$}\\
          X(t) & \text{if \ $n/2<t< n/2+1,$}\\
          X(n/2+1) & \text{if \ $t\geq n/2+1,$}\\
        \end{cases}
        \qquad n\in\NN\cup\{0\},
 \]
 and hence for all $n\in\NN\cup\{0\},$ the process $\{X^{(n)}(t) : t\geq 0\}$ has a continuous modification
 which is locally $\gamma$-H\"older continuous for every $\gamma\in(0,q/p).$
The desired property of $\{X(t) : t\geq 0\}$ follows by that
 $[0,\infty) = \bigcup_{n\in\NN\cup\{0\}}[n/2,n/2+1].$

Now, let $p$ be a positive even number.
Using the stationary increments property
 and the relation between moments and cumulants we obtain
 \begin{equation}\label{Kolmo1}
  \EE\big\vert X(t)-X(s)\big\vert^p
  =\EE  (X(t-s))^p
  =\sum_\Pi\prod_{B\in\Pi}\cum[n_B](X(t-s)),
 \end{equation}
 where $\Pi$ runs through the list of all partitions of a set of size \!$p,$\!
 $B\in\Pi$ means that \!$B$\! is one of the blocks into which the set is partitioned
 given the partition \!$\Pi,$
 $n_B$ is the size of the set \!$B$\! (in notation: $n_B=\vert B\vert$\!)
 and $\cum[k](X(t-s))$ denotes the cumulant of order \!$k$\! of  $X(t-s),$
 see, e.g., Shiryaev \cite[page 292, formula (46)]{Shi}.
Observe that $\EE X(t-s)=0$ implies that for any \!$\Pi$\!
and any one element block \ $B\in\Pi,$ \ we have $\cum[n_B](X(t-s))=\cum[1](X(t-s))=\EE X(t-s)=0,$
 which yields that for any partition $\Pi$ the sum on the right hand side of \eqref{Kolmo1}
 has at most $p/2$ terms different from $0.$
Using the dilative stability relation (\ref{dsrel}) we can continue (\ref{Kolmo1}) as follows:
\begin{align}
\EE\big\vert X(t)-X(s)\big\vert^p
&=\sum_\Pi\prod_{B\in\Pi}\cum[n_B](X(t-s))\notag\\[5pt]
&=\sum_\Pi\prod_{B\in\Pi} (t-s)^{(H-\delta/2)n_B+\delta}\cum[n_B](X(1))\notag\\[5pt]
&= (t-s)^{(H-\delta/2)p}
 \sum_\Pi\bigg( (t-s)^{\delta\vert\Pi\vert}\prod_{B\in\Pi}\cum[n_B](X(1))\bigg).
\label{parti}
\end{align}
It is important to observe that for each $\Pi,$ the product $\prod_{B\in\Pi}\cum[n_B](X(1))$
 is nonnegative.
Indeed, if $n_B$ is even, then $\cum[n_B](X(1))\geq0,$ since the cumulant of order
greater than or equal to 2 of $X(1)$ is the moment of the same order of the L\'evy measure in the
 L\'evy--Khinthine representation of the distribution of $X(1)$ (plus the variance of the Gaussian component
 if $n_B=2$), see Steutel and Van Harn \cite[Chapter IV, Theorem 7.4]{Ste}.
A cumulant $\cum[n_B](X(1))$ can be negative for an odd $n_B,$ however, the number of blocks $B\in\Pi,$ with an odd size $n_B,$ must be even, since \!$p$\! is even.
This yields the nonnegativity of $\prod_{B\in\Pi}\cum[n_B](X(1)).$
At this point the proof separates into three cases corresponding to the three parts of the statement of the theorem.

$\bullet$ If $\delta<0,$ then (\ref{parti}) can be continued in the following way:
\begin{align*}
\EE\big\vert X(t)-X(s)\big\vert^p
&= (t-s)^{(H-\delta/2)p}
   \sum_\Pi\bigg((t-s)^{\delta\vert\Pi\vert}\prod_{B\in\Pi}\cum[n_B](X(1))\bigg)\\[5pt]
&\leq (t-s)^{(H-\delta/2)p}
 \sum_\Pi\Big( (t-s)^{\delta p/2}\prod_{B\in\Pi}\cum[n_B](X(1))\Big)\\[5pt]
&=\EE  (X(1))^p (t-s)^{Hp},
\end{align*}
where at the inequality we used the facts that $0\leq t-s<1$
and $\vert\Pi\vert\leq p/2$ for all $\Pi$ not having one-element blocks; and the last equality
 follows by \eqref{Kolmo1}.
Choosing $p>1/H,$ we conclude, by the Kolmogorov--Chentsov theorem, that
$\{X(t),\,t\geq0\}$ has a continuous modification which is locally \!$\gamma$\!-H\"older
 continuous for every
\[
  0<\gamma<\frac{Hp-1}{p}=H-\frac1p\,.
\]
By letting $p\to\infty$ we have finished the proof of the case (i).

$\bullet$ If $0\leq\delta<2H,$ the proof proceeds similarly.
Using (\ref{parti}) we obtain
\begin{align*}
\EE\big\vert X(t)-X(s)\big\vert^p
&= (t-s)^{(H-\delta/2)p}
    \sum_\Pi\bigg((t-s)^{\delta\vert\Pi\vert}\prod_{B\in\Pi}\cum[n_B](X(1))\bigg)\\[5pt]
&\leq (t-s)^{(H-\delta/2)p}\sum_\Pi\Big( (t-s)^0\prod_{B\in\Pi}\cum[n_B](X(1))\Big)\\[5pt]
&=\EE (X(1))^p (t-s)^{(H-\delta/2)p}.
\end{align*}
For $p>1/(H-\delta/2),$ the Kolmogorov--Chentsov theorem ensures that
$\{X(t),\,t\geq0\}$ has a continuous modification which  is locally \!$\gamma$\!-H\"older
 continuous for every
\[
    0<\gamma<\frac{(H-\delta/2)p-1}{p}=H-\frac\delta2-\frac1p\,.
\]
Letting $p\to\infty$ as above, we obtain the statement (ii).

$\bullet$ If $\delta=2H$ and $H>1/2,$ it is enough to use the second moment:
\[
\EE\big\vert X(t)-X(s)\big\vert^2
   = \EE(X(t-s))^2 = \EE (X^{\circledast (t-s)^{2H}}(1))^2
   =\EE (X(1))^2 (t-s)^{2H},
\]
 to conclude that (by the Kolmogorov--Chentsov theorem) $\{X(t), t\geq 0\}$ has a
continuous modification which  is locally \!$\gamma$\!-H\"older continuous for every
\[
   0<\gamma<\frac{2H-1}2=H-\frac12,
\]
 which is the statement (iii).

Finally, the statements $(\widetilde i)$--$(\widetilde{iii})$ follow from  part (ii) of Lemma \ref{Lem2},
 using also that for a $(H,\delta)$-dilatively stable process with stationary increments,
 the range of the parameter \!$H$\! is  $(0,1]$ (see the Introduction).
\proofend

\begin{Rem}
The key in the proof of Theorem \ref{Holder-thm} was the Kolmogorov condition,
which is known to be not a necessary condition for
having a continuous modification of a stochastic process, see Stoyanov \cite[p. 219]{Sto}.
It is not a necessary condition for  having a locally H\"older continuous modification either,
 as one can see using a modification of the counterexample of Stoyanov \cite[p. 220]{Sto}.
Namely, if $\{W(t),t\geq0\}$ is a FBM with parameter \!$H,$\! then
 $X(t)\doteq\exp(W^3(t)),$  $t\geq 0,$ has infinite moments,
 hence the Kolmogorov condition does not make sense.
However, as $\{W(t),t\geq0\}$  has a continuous modification which is
 locally \!$\gamma$\!-H\"older continuous for every $\gamma\in(0,H)$
 (see  Remark \ref{REM_ss1}), we have
  for every bounded subinterval $J\subseteq[0,\infty),$
\[
\vert X(t)-X(s)\vert
\leq c_1\big\vert W^3(t)-W^3(s)\big\vert
\leq c_2\vert W(t)-W(s)\vert
\leq c_3\vert t-s\vert^\gamma,
  \quad s,t\in J,
\]
with some (random) constants $c_1,c_2,c_3,$
hence  $\{X(t),t\geq0\}$ is a.s. locally  \!$\gamma$\!-H\"older continuous for every $\gamma\in(0,H).$
 (Here one can choose universal constants $c_i,$ $i=1,2,3,$ i.e., which do not depend on the specific
 choices of $s,t\in J,$ since $W$ is almost surely bounded on the bounded interval $J.$\!)
Therefore it is reasonable to ask how strong the assertions  $(\widetilde {ii})$ and $(\widetilde{iii})$
 of Theorem \ref{Holder-thm} are.
By  Benassi et al. \cite[Proposition 3.2]{Ben}, if $\gamma>H-1/2,$ then  on any interval
 the sample paths of a non-Gussian FLP with parameter \!$H$\!
 are not \!$\gamma$\!-H\"older continuous
 with positive probability $p>0.$
Furthermore, if the L\'evy measure (control measure) of the L\'evy process
in the defining integral of the non-Gussian FLP is infinite, then $p=1.$
This example shows that part  $(\widetilde {ii})$ of Theorem \ref{Holder-thm}
 cannot be strengthened in general, in the sense that, as the above shows, there exists
 an $(H,\delta)$\!-dilatively stable process with stationary increments and zero mean for which
 $0\leq \delta<2H$ and $\Hexp_X=H-\delta/2$  (namely, a non-Gaussian FLP with parameter \!$H$\!
 for which $\delta=1$), i.e., the left endpoint of the interval
 $[H-\delta/2,H]$ can be reached.
The authors do not know whether the right endpoint of the interval $[H-\delta/2,H]$ can be reached.
We have only a partial result in case of $\delta=0.$
Namely, by part $(\widetilde{ii})$ of Theorem \ref{Holder-thm} with $\delta=0,$ we
 see that if \!$X$\! is a \!$H$\!-self-similar process with stationary increments,
 having non-Gaussian, infinitely divisible finite-dimensional distributions,
finite moments of all orders, and a Gaussian component,
then it a.s. admits a local H\"older exponent \!$H.$
\end{Rem}

Next we formulate two corollaries of Theorem \ref{Holder-thm}.

\begin{Cor}
For the LISOU, LISCBI and  LISDLG processes with parameter \!$H$\!
(see Example \ref{LISexample})  one can apply parts (i) and $(\widetilde {i})$ of
 Theorem \ref{Holder-thm} (remember, $1/2<H<1,$ hence $\delta=2H-2<0$\!),
 hence the local H\"older exponent  $\Hexp_X$ of these processes
 is a.s. \!$H.$\!
On the other hand, the  non-Gaussian FLP with parameter \!$H$\!
 (but possibly with a Gaussian component, see again Example \ref{LISexample})
 is a.s. locally \!$\gamma$\!-H\"older continuous
 for every $0<\gamma<H-1/2$ by part (ii) of Theorem \ref{Holder-thm}
 (since we have $1/2<H<1,$ hence $\delta=1<2H$\!),
 which is known by Benassi et al. \cite[Proposition 3.2]{Ben} or Marquardt \cite[Theorem 4.3]{Mar}.
Particularly, all four processes have continuous modifications.
\end{Cor}

\begin{Cor} \label{C-beliseg}
Let $\{X(t),\,t\geq0\}$ be an $(H,\delta)$\!-dilatively stable process
 with stationary increments and zero mean such that either $\delta{\,\neq\,}2H$ or $H{\,>\,}1/2.$
Then  $\{X(t),\,t\geq0\}$ has  a continuous modification.
\end{Cor}

\begin{Rem}\label{REM_ss1}
The self-similar analogue of Theorem \ref{Holder-thm} is a well-known
 and easy consequence of the Kolmogorov--Chentsov theorem.
Namely, for a $H$\!-self-similar process $\{X(t),\,t\geq0\}$ with stationary
 increments and finite moments of all orders,
 we have $\EE\vert X(t)-X(s)\vert^p=\vert t-s\vert^{pH}\EE\vert X(1)\vert^p,$ $p>0,$
 from which we obtain, by the same way as in the proof of  Theorem \ref{Holder-thm}, that
 $\{X(t),\,t\geq0\}$  has  a continuous modification, the sample paths of which
 are locally \!$\gamma$\!-H\"older continuous for every  $\gamma\in(0,H).$
For completeness, we note that the zero mean condition in Theorem \ref{Holder-thm}
 is automatically satisfied by a $H$-self-similar process with $H\ne 1,$
 stationary increments and finite  absolute moment,
 see, e.g., Vervaat \cite[Auxiliary Theorem 3.1]{Ver}.
In the special case $H=1$ (and $\delta=0$\!) it is not sure that we have a zero mean process,
 but it is a degenerate case, i.e., $X(t)=tX(1),$ $t\geq 0,$ a.s., see, e.g., Igl\'oi
 \cite[Theorem 2.7.1 (SS) 2)]{Igl}, and hence the corresponding assertions of parts
 $(ii)$ and $(\widetilde{ii})$ of Theorem \ref{Holder-thm} hold readily.

Further, assuming that $\PP(X(1)=0)=0$, i.e., the distribution of $X(1)$ has no atom at zero,
 the proof of Lemma \ref{Lem1} in the case of $\delta=0$ shows that the local H\"older exponent
 $\Hexp_X$\! equals a.s. \!$H.$
\end{Rem}

\section{Singularity of the distributions}
\label{singulsect}

In what follows  $C(I)$ will denote the set of continuous functions
on  a closed interval $I\subseteq[0,\infty)$ with the local
uniform topology and Borel $\sigma$\!\!-algebra  $\cB(C(I)).$
If  the processes $\{X_1(t),\,t\in I\}$ and $\{X_2(t),\,t\in I\}$ have
 sample paths  in $C(I)$ a.s.,
then we say that they are singular on \!$I$\! (in notation: $X_1\,\bot\, X_2$\!),
if their distributions $\P_{X_1}$ and $\P_{X_2}$ on $(C(I),\cB(C(I)))$ are singular
 (in notation:  $\P_{X_1}\,\bot\, \P_{X_2}$\!), i.e., there exists a set $A\in\cB(C(I))$
 such that $\P_{X_1}(A)=1$ and $\P_{X_2}(A)=0.$

\begin{Thm}
\label{Thm1}
Let \!$I\subseteq [0,\infty)$\! be a closed interval
and $\{X_1(t),\,t\geq0\},$ $\{X_2(t),\,t\geq0\}$  be dilatively stable processes with parameters
$(\alpha_1,\delta_1)$ and $(\alpha_2,\delta_2),$  respectively, both processes with zero mean,
and having sample paths  in $C(I)$ a.s.
Assume  that one of the following two conditions holds:\newline
$\bullet$  $\inf\{t : t\in I\}=0$ (i.e., the left endpoint of \!$I$\! is zero) and the above
 two processes satisfy Assumption \ref{Assum_delta+}.\newline
$\bullet$  $\sup\{t : t\in I\}=\infty$ (i.e., $I$ is unbounded) and the above
 two processes satisfy Assumption \ref{Assum_delta-}.\newline
Then $\alpha_1\neq\alpha_2$ implies  $X_1\,\bot\, X_2.$
\end{Thm}

\noindent{\bf Proof.} Assume that the left endpoint of \!$I$\!  is zero,
 and Assumption \ref{Assum_delta+} holds  for the processes $X_1$ and $X_2.$
Let $(t_n)_{n\in\mathbb{N}}$ be a sequence for which the assertions of Lemma \ref{Lem1} hold,
 and define the following two subsets of $C(I)$\!:
\begin{align*}
 A_i\doteq\bigg\{f\in C(I):
 &\limsup_{n\to\infty}\frac{\vert f(t_n)\vert}{t_n^\kappa}=0\,\
 \text{for all}\,\ \kappa\in\mathbb{Q}\cap(0,\alpha_i)\\
 &\;\;\text{and}\,\
 \limsup_{n\to\infty}\frac{\vert f(t_n)\vert}{t_n^\kappa}=\infty\,\
 \text{for all}\,\ \kappa\in\mathbb{Q}\cap(\alpha_i,\infty)\bigg\},\quad i=1,2,
\end{align*}
 where $\QQ$  denotes the set of rational numbers.
Considering the decompositions
 \begin{align*}
   A_i = \bigcap_{\kappa\in\QQ\cap(0,\alpha_i)}
            &\left(\bigcap_{m=1}^\infty \bigcup_{p=1}^\infty \bigcap_{r=p}^\infty
                 \Big\{ f\in C(I) : \frac{\vert f(t_r)\vert}{t_r^\kappa}<\frac{1}{m}\Big\}
             \right)\\
            &\bigcap
           \bigcap_{\kappa\in\QQ\cap(\alpha_i,\infty)}
             \left(
                \bigcap_{m=1}^\infty
              \bigcap_{p=1}^\infty\bigcup_{r=p}^\infty\Big\{f\in C(I):\frac{\vert f(t_r)\vert}{t_r^\kappa} >m\Big\}
            \right), \qquad i=1,2,
 \end{align*}
 and using that for each $\kappa\in\QQ$ and $n\in\NN,$
 the mapping $C(I)\ni f\mapsto f(t_n)/t_n^\kappa$ is continuous, we get $A_i\in\cB(C(I)),$ $i=1,2.$
One can argue in another way, namely, by Lemma \ref{Lem1},
 for all $\kappa\in\QQ\cap(0,\alpha_i)$ [resp. $\kappa\in\QQ\cap(\alpha_i,\infty)$],
 we have
 \[
   \Big\{ f\in C(I) : \limsup_{n\to\infty}\frac{\vert f(t_n)\vert}{t_n^\kappa}=0 \Big\}
   \qquad
   \Big[ \;\text{resp.} \quad \Big\{ f\in C(I) :
              \limsup_{n\to\infty}\frac{\vert f(t_n)\vert}{t_n^\kappa}=\infty \Big\}\; \Big]
 \]
 is in $\cB(C(I)).$
Hence, by Lemma \ref{Lem1}, $\P_{X_i}(A_i)=1,$ $i=1,2.$
Since $\alpha_1\neq\alpha_2$ we have $A_1\cap A_2=\emptyset,$  hence
 the assertion follows.

In the other case the proof is analogous, but we have to refer
to Lemma \ref{Leminfty} instead of Lemma \ref{Lem1}.
\proofend

\medskip

The next theorem is the counterpart of Theorem \ref{Thm1} for processes with stationary increments,
 in which case the closed interval \!$I$\! can be arbitrary.

\begin{Thm}
\label{Thm2}
Let $I\subseteq[0,\infty)$ be a closed  interval,
and $\{X_1(t),\,t\geq0\},$ $\{X_2(t),\,t\,\geq0\}$  be dilatively stable processes
 with stationary increments and parameters $(H_1,\delta_1)$  and $(H_2,\delta_2),$ respectively,
 both processes with zero mean and having sample paths in $C(I)$ a.s. such that
 they satisfy Assumption \ref{Assum_delta+}.
Then  $H_1\neq H_2$ implies $X_1\,\bot\, X_2.$
\end{Thm}

\noindent{\bf Proof.}
One can argue similarly to the proof of Theorem \ref{Thm1}.
Namely, let $t_0\in I,$\!
$(t_n)_{n\in\mathbb{N}}$ be  a sequence for which the assertions of Lemma \ref{Lem2} hold,
 and for $i=1,2,$
\begin{align*}
 A_i\doteq\bigg\{f\in C(I):
 &\limsup_{n\to\infty}\frac{\vert f(t_n) -  f(t_0) \vert}{\vert t_n - t_0\vert^\kappa}=0\,\
 \text{for all}\,\ \kappa\in\mathbb{Q}\cap (0,\alpha_i)\\
 &\;\;\text{and}\,\
 \limsup_{n\to\infty}\frac{\vert f(t_n) - f(t_0) \vert}{\vert t_n - t_0\vert^\kappa}=\infty\,\
 \text{for all}\,\ \kappa\in\mathbb{Q}\cap(\alpha_i,\infty)\bigg\}.
\end{align*}
Hence the assertion follows by Lemma \ref{Lem2}.
\proofend

\begin{Rem}
If either $\delta\ne 2H$ or $H>1/2,$ then in Theorem \ref{Thm2}
the assumption of having a.s. continuous sample paths is automatically fulfilled.
Indeed, by Corollary \ref{C-beliseg}, apart from the exceptional case $\delta{\,=\,}2H$
and $0{\,<\,}H{\,\leq\,}1/2,$ every dilatively stable process with stationary increments
and zero mean has a continuous modification.
\end{Rem}

\begin{Ex}
By Theorem \ref{Thm2}, any two LISOU processes with different parameters \!$H$\! are singular on any closed
 interval $I\subseteq[0,\infty),$ because Assumption \ref{Assum_delta+} is trivially fulfilled
 (since $\delta=2H-2<0$\!).
The same is true for the LISCBI, and particularly, for the LISDLG processes,
 since they have the same parameter of dilative stability as the LISOU process.
It also follows that any two of these processes of three types are singular if their parameters
  \!$H$\! are different.
However, for the non-Gaussian FLP we have $\delta{\,=\,}1{\,\geq\,}0,$ hence
 Assumption  \ref{Assum_delta+} is non-trivial, and Theorem \ref{Thm2} states that
 two non-Gaussian FLPs, both having Gaussian components, are singular
 if their parameters \!$H$\! are different.
\end{Ex}

Theorem \ref{Thm2} applies to the case $\delta{\,=\,}0,$ i.e., for $(\alpha,0)$\!-dilatively
stable processes too. These processes are exactly those self-similar, infinitely divisible processes,
which have finite moments of all orders and non-Gaussian one-dimensional distributions (except for that at zero).
But of course, there exist self-similar processes without these properties, e.g. the FBM.
For this reason, Theorem \ref{Thm2} does not automatically apply to self-similar processes in general.
However, let us observe that not all of the three properties above are utilized in the proof:
neither the non-Gaussianity, nor the moments (or cumulants) of order higher than two are in use
(recall the  proofs of Lemma \ref{Lem1}  and Theorem \ref{Thm2} in the case of $\delta=0$\!),
 and in case of $\delta=0$ we had nothing to do with infinite divisibility
(since the convolution exponent $T^\delta$ in \eqref{dsrel} is 1 in this case).
In fact, not even finite second moments are needed
for Theorem \ref{Thm2} to apply to the self-similar case.
Indeed, the only place where moments appear  in case of $\delta=0,$
 is the Markov inequality in the proof of part (i) of Lemma \ref{Lem1},
 where one can use, e.g., the absolute moment instead of the second order one.
Let us also observe that when $\delta=0,$ Assumption \ref{Assum_delta+}
 (i.e., the existence of a Gaussian component)  is utilized only in a way that
 it follows that the  dilatively stable process in question does not have an atom at zero.
So, if we replace dilatively stable processes by self-similar ones possessing the following properties,
 then Theorem \ref{Thm2} remains true, and reads as follows.

\begin{Thm}
\label{Thm2ss}
Let $I\subseteq[0,\infty)$ be a closed interval,
 and $\{X_1(t),\,t\geq0\},$ $\{X_2(t),\,t\,\geq0\}$ be self-similar processes
 with stationary increments and  parameters $H_1,H_2\in(0,1),$
 both processes with finite absolute moment and having sample paths in $C(I)$  a.s.
 such that neither the distribution of $X_1(1)$ nor that of $X_2(1)$ has an atom at zero.
Then  $H_1\neq H_2$ implies $X_1\,\bot\, X_2.$
\end{Thm}

Note that under the conditions Theorem \ref{Thm2ss}, the processes $X_1$ and $X_2$
 have zero mean, see, Vervaat \cite[Auxiliary Theorem 3.1]{Ver}.

The most important particular case of Theorem \ref{Thm2ss}  sounds as follows.

\begin{Cor}
Two FBMs with different parameters \!$H$\! are singular.
\end{Cor}

In fact, this latter result is known, see Prakasa Rao \cite{Rao},
 where the proof is based on a Baxter type theorem of Kurchenko \cite{Kur}.

\begin{Rem}
A consequence of the above theorems on singularity
(which are in fact based on Lemmas \ref{Lem1} [or \ref{Leminfty}] and \ref{Lem2})
 is that the parameter \!$\alpha,$\! or \!$H$\! in the stationary increments case,
 can be estimated without error, as long as we have a continuous-time sample path available
 (or at least its values at the time points of a sequence $(t_n)_{n\in\NN}$
 appearing in Lemmas \ref{Lem1} [\ref{Leminfty}] and \ref{Lem2}).
 Note also that in Lemmas \ref{Lem1} [\ref{Leminfty}] and \ref{Lem2} not just the existence
 of an appropriate sequence $(t_n)_{n\in\NN}$ is guaranteed, but we also give an example
 for such a sequence, which is important from the point of view of practical applications.
\end{Rem}

Finally, we recall a known result about variation of sample paths of self-similar processes,
 the generalization of which may serve as a future task.

\begin{Rem}
If $\{X(t), t\geq 0\}$ is a \!$H$\!-self-similar process with stationary increments,
 finite  absolute moment and $H<1$ then the sample paths of \!$X$\! have no bounded
 variation on any (bounded) interval a.s., see Vervaat \cite[Theorem 3.3]{Ver}.
As a possible future task one can study variation of sample paths of dilatively stable processes.
\end{Rem}

\appendix

\section{On the definition of dilatively stability}\label{appendix}

First we note that, at the first view, the Definition \ref{DefDS} of dilative stability is a little bit
 different from Definition 2.1.3 in Igl\'oi \cite{Igl}, since the right continuity of the
$n$\!-th order (\!$n\geq 2,$ $n\in\NN$\!) cumulant function $c_n:[0,\infty)\to\RR,$
 $c_n(t):=\cum[n](X(t)),$ $t\geq 0,$ is not supposed.
However, it follows since $c_n(t)=t^{(\alpha-\delta/2)n+\delta}c_n(1),$ $t\geq 0.$
Even so, this does not mean that Definition 2.1.3 in Igl\'oi \cite{Igl} contains a superfluous condition
 (since it defines dilative stability in a more general setting) and it turns out to be equivalent to our definition,
 see Igl\'oi \cite[Theorem 2.2.1]{Igl}.

Next we note that there is a slight redundancy in Definition \ref{DefDS} in the sense that
 the property that a dilatively stable process starts from $0$ follows from the other properties.
More precisely, if a process $\{X(t),\,t\geq0\}$ satisfies the properties listed in Definition \ref{DefDS}
 except that it starts from $0,$ then $\P(X(0)=0)=1.$
Indeed, by \eqref{dsrel}, $X(0)\sim T^{\alpha-\delta/2} X^{\circledast T^\delta}(0)$
 (where $\sim$ denotes equality in distribution), which yields that
 \begin{align*}
   \EE X(0) = T^{\alpha-\delta/2}\EE X^{\circledast T^\delta}(0)
             = T^{\alpha-\delta/2} T^{\delta}\EE X(0)
             = T^{\alpha+\delta/2}\EE X(0),\qquad T>0,
 \end{align*}
 and
 \begin{align*}
   \D^2X(0) = T^{2\alpha-\delta}\D^2X^{\circledast T^\delta}(0)
             = T^{2\alpha-\delta} T^{\delta}\D^2X(0)
             = T^{2\alpha}\D^2X(0), \qquad T>0.
 \end{align*}
If $\delta\ne-2\alpha,$ then this implies that $\EE X(0)=0$ and $\D^2X(0)=0,$ yielding that
 $\P(X(0)=0)=1.$
If $\delta=-2\alpha,$ then, by \eqref{dsrel}, $X(0)\sim T^{2\alpha}X^{\circledast T^{-2\alpha}}(0)$ for all $T>0,$
 or equivalently $cX(0)\sim X^{\circledast c}(0)$ for all $c>0,$
which yields that the distribution of $X(0)$ is strictly $1$\!-stable.
Using that a non-degenerate (strictly) $1$\!-stable distribution does not have a finite first moment,
 our assumption that $X(0)$ has (finite) moments of all orders implies that $\PP(X(0)=C)=1$ with some $C\in\RR.$
For completeness, we also note that there is a slight redundancy in Definition \ref{DefSS} too,
 in the sense that if a process $\{X(t),\,t\geq0\}$ satisfies the scaling property \eqref{ssrel},
 then $\P(X(0)=0)=1.$

Finally, we call the attention that the
non-Gaussianity condition in Definition \ref{DefDS}, namely, that $X(1)$
(or, equivalently, $X(t),$ for some $t>0$)
is non-Gaussian
 is crucial in the sense that it is extensively used in the proofs
 both in the present paper and in Igl\'oi \cite{Igl}.
Note also that this condition ensures that all the higher-dimensional distributions are also non-Gaussian.

\section*{Acknowledgements}
M. Barczy has been supported by the NKTH-OTKA-EU FP7 (Marie Curie action)
 co-funded 'MOBILITY' Grant No. OMFB-00610/2010, and by the Hungarian Scientific Research
 Fund under Grant No.\ OTKA T--079128.

\end{document}